\documentclass{amsproc}

\usepackage{csquotes}
\usepackage{xspace}
\usepackage{pgf,tikz}
\usepackage{booktabs}
\usepackage{amsmath,amssymb,amsthm}
\usepackage{mathbbol}
\usepackage{url}

\usetikzlibrary{calc,3d,shapes,fit,arrows}
\usetikzlibrary{decorations.markings}
\usetikzlibrary{patterns}
\usetikzlibrary{intersections}

\globalcolorstrue

\definecolor{cobaltblue}{cmyk}{100,54,0,0}
\colorlet{blue}{cobaltblue}
\definecolor{strawberry}{RGB}{197,70,68}
\colorlet{red}{strawberry}
\colorlet{myred}{strawberry}

\tikzset{VertexStyle/.style = {
    circle, draw,
    text           = black,
    inner sep      = 2pt,
    outer sep      = 0pt,
    minimum size   = 20 pt}}

\tikzset{BlackDot/.style = {
    VertexStyle,
    fill = black}}
\tikzset{GrayDot/.style = {
    VertexStyle,
    fill = gray}}
\tikzset{WhiteDot/.style = {
    VertexStyle,
    fill = white}}
\tikzset{RedDot/.style = {
    VertexStyle,
    fill = red}}
\tikzset{GreenDot/.style = {
    VertexStyle,
    fill = green}}
\tikzset{BlueDot/.style = {
    VertexStyle,
    fill = blue}}

\tikzset{EdgeStyle/.style = {
    color = black!70,
    thick,
    >=latex,
  }}

\tikzset{MatchingEdgeStyle/.style = {
    color = red,
    ultra thick,
    >=latex,
  }}

\tikzset{ArcStyle/.style = {
    EdgeStyle,
    ->
  }}

\tikzset{LabStyle/.style = {
    fill=white,
  }}

\tikzset{LoopStyle/.style = {text = red!80}}

\tikzset{KoordStyle/.style = {->, >=latex, color=black!40, thick}}
\tikzset{GitterStyle/.style = {very thin,color=black!40,join=round}}

\tikzset{FillPoly/.style = {fill opacity = 0.4, fill=strawberry!70}}
\tikzset{DrawPoly/.style = {very thick, draw=strawberry}}
\tikzset{Poly/.style = {FillPoly, DrawPoly}}
\tikzset{Hyperplane/.style = {ultra thick, draw=gray}}
\tikzset{HyperplaneLabel/.style = {circle, draw=black, fill = white}}
\tikzset{Unimodular/.style = {DrawPoly, fill opacity = 0.3, fill=red!90}}

\tikzstyle{pseudovertex} = [circle, scale=0.6, fill=strawberry]
\tikzstyle{generator} = [circle, scale=0.6, draw=strawberry, fill=white]
\tikzstyle{vertex} = [rectangle, scale=0.8, draw=strawberry, fill=white]

\newcommand{\StandardGitter}{
  \draw[GitterStyle] (-2,-2) -- (2,-2) -- (2,2) -- (-2,2) -- cycle;
  \foreach \i in {-1,0,1}{
    \draw[GitterStyle](\i,-2) -- (\i,-1.8);
    \draw[GitterStyle](\i,2) -- (\i,1.8);
    \draw[GitterStyle](-2,\i) -- (-1.8,\i);
    \draw[GitterStyle](2,\i) -- (1.8,\i);
    \small
    \node at (-2,\i)[left] {$\i$};
    \node at (\i,-2)[below] {$\i$};
  }
  \foreach \i in {-1.5, -0.5, 0.5, 1.5}{
    \draw[GitterStyle](\i,-2) -- (\i,-1.9);
    \draw[GitterStyle](\i,2) -- (\i,1.9);
    \draw[GitterStyle](-2,\i) -- (-1.9,\i);
    \draw[GitterStyle](2,\i) -- (1.9,\i);
  }
}
  
\newcommand{\PolygonGitter}{
  \draw[GitterStyle] (-3,-2) -- (6,-2) -- (6,3) -- (-3,3) -- cycle;
  \foreach \i in {-2,-1,...,5}{
    \draw[GitterStyle](\i,-2) -- (\i,-1.8);
    \draw[GitterStyle](\i,3) -- (\i,2.8);
    \node at (\i,-2)[below] {$\i$};
  }
  \foreach \i in {-1,0,1,2}{
    \draw[GitterStyle](-3,\i) -- (-2.8,\i);
    \draw[GitterStyle](6,\i) -- (5.8,\i);
    \small
    \node at (-3,\i)[left] {$\i$};

  }
  \foreach \i in {-2.5, -1.5, ..., 5.5}{
    \draw[GitterStyle](\i,-2) -- (\i,-1.9);
    \draw[GitterStyle](\i,3) -- (\i,2.9);
  }
  \foreach \i in {-1.5, -0.5, ..., 2.5}{
    \draw[GitterStyle](-3,\i) -- (-2.9,\i);
    \draw[GitterStyle](6,\i) -- (5.9,\i);
  }
}

\colorlet{verylightgray}{black!10!}
\colorlet{lightgray}{black!25!}
\colorlet{myblue}{blue!75!}
\tikzstyle{blackdot} = [circle,draw,fill=black,scale=0.5]
\tikzstyle{whitedot} = [circle,draw,fill=white,scale=0.5]
\tikzstyle{edge} = [draw,-,black]
\tikzstyle{mediumedge} = [draw,thick,-,black,join=round]
\tikzstyle{thickedge} = [draw,ultra thick,-,black,join=round]

\newcommand\DrawMaxTropLine[4]{
  \draw[#1] (#2,#3) -- (#2-#4,#3); 
  \draw[#1] (#2,#3) -- (#2,#3-#4); 
  \draw[#1] (#2,#3) -- (#2+#4,#3+#4);
}

\newtheorem{theorem}{Theorem}
\newtheorem{proposition}[theorem]{Proposition}

\newtheorem{lemma}[theorem]{Lemma}
\theoremstyle{definition}

\newtheorem{example}[theorem]{Example}
\theoremstyle{remark}
\newtheorem{remark}[theorem]{Remark}

\newtheorem{observation}[theorem]{Observation}

\numberwithin{equation}{section}

\newcommand\KK{{\mathbb K}}
\newcommand\LL{{\mathbb L}}
\newcommand\NN{{\mathbb N}}
\newcommand\QQ{{\mathbb Q}}
\newcommand\RR{{\mathbb R}}
\newcommand\TT{{\mathbb T}}
\newcommand\ZZ{{\mathbb Z}}
\newcommand\cL{{\mathcal L}}

\newcommand\cS{{\mathcal S}}
\newcommand\cT{{\mathcal T}}
\newcommand\1{{\mathbb 1}}
\newcommand\TP{{\mathbb{TP}}}
\newcommand\TTmax{\TT_{\max}} 

\newcommand\SetOf[2]{\left\{#1 \mid #2\right\}}
\newcommand\smallSetOf[2]{\{#1 \mid #2\}}

\newcommand\torus[1]{\RR^{#1}/\RR\1} 
\DeclareMathOperator\trop{trop} 
\DeclareMathOperator\tconv{tconv} 
\DeclareMathOperator\supp{supp} 
\DeclareMathOperator\TLP{TLP} 
\DeclareMathOperator\Eig{Eig}
\newcommand\minodot{\mathbin{\underline{\odot}}}
\newcommand\maxodot{\mathbin{\overline{\odot}}}
\newcommand\minoplus{\mathbin{\underline{\oplus}}}
\newcommand\maxoplus{\mathbin{\overline{\oplus}}}
\newcommand\val{\operatorname{ord}} 

\newcommand\puiseuxseries[2]{{#1}\{\hskip-.25em\{#2\}\hskip-.25em\}}

\newcommand\scp[2]{\langle #1,#2 \rangle} 
\newcommand\tscp[2]{\scp{#1}{#2}_{\rm trop}} 

\newcommand\NP{\text{\rm NP}}
\newcommand\coNP{\text{\rm co-NP}}

\newcounter{tikzbrace}
\newcommand{\tikzmark}[1]{\tikz[overlay,remember picture] {\node (brace-\thetikzbrace) {};}\stepcounter{tikzbrace}}
\newcommand{\insertbigbrace}[1]{%
\begin{tikzpicture}[remember picture, overlay]
\draw[thick] 
	let \n1 = {\thetikzbrace - 2},
		\n2 = {\thetikzbrace - 1},   
		\p1 = (brace-\n1),
		\p2 = (brace-\n2),
		\n3 = {max(\x1,\x2)},
		\p3 = ($(\n3,\y1) + (0.2,0.4)$),
		\p4 = ($(\n3,\y2) + (0.2,-0.2)$),
		\p5 = (-0.1,0) in 
	(\p3) ++ (\p5) -- (\p3) -- node[right=1ex]{#1} (\p4) -- ++ (\p5);
\end{tikzpicture}}

\DeclareMathOperator\SL{SL}
\DeclareMathOperator\Gr{Gr}
\DeclareMathOperator\Dr{Dr}
\DeclareMathOperator\TGr{TGr}

\newcommand\mptopcom{\texttt{mptopcom}\xspace}
\newcommand\TOPCOM{\texttt{TOPCOM}\xspace}

\newcommand\Maple{\texttt{Maple}\xspace}

\newcommand\doi[1]{\href{http://dx.doi.org/#1}{\texttt{doi:#1}}}

\begin{document}

\title{Developments in tropical convexity}

\author{Michael Joswig} 

\address[Michael Joswig]{
  Technische Universität Berlin,
  Chair of Discrete Mathe\-ma\-tics/Geo\-me\-try \\
  and Max-Planck Institute for Mathematics in the Sciences, Leipzig
}
\email{joswig@math.tu-berlin.de}
\thanks{%
  Supported by Deutsche Forschungsgemeinschaft (EXC 2046: \enquote{MATH$^+$} and SFB-TRR 195: \enquote{Symbolic Tools in Mathematics and their Application}).
  The hospitality of the Hausdorff Institute for Mathematics in Bonn, during the Trimester Program \enquote{Discrete Optimization} is gratefully acknowledged.
}

\date{\today}

\dedicatory{Dedicated to Bernd Sturmfels on the occasion of his 60th birthday.}


\begin{abstract}
  The term \enquote{tropical convexity} was coined by Develin and Sturmfels who published a landmark paper with that title in 2004.
  However, the topic has much older roots and is deeply connected to linear and combinatorial optimization and other areas of mathematics.
  The purpose of this survey is to sketch how that article contributed to shaping the field of tropical geometry as we know it today.
\end{abstract}

\maketitle

\section{Introduction} 

In their article \cite{SpeyerSturmfels:2009} Speyer and Sturmfels wrote \enquote{$\ldots$ the tropical approach in mathematics [is] now an integral part of geometric combinatorics and algebraic geometry.  It has also expanded into mathematical physics, number theory, symplectic geometry, computational biology, and beyond.}
Why is this the case, and how did that happen?
I think the reason is that the principles underlying tropical geometry are at the same time elementary and nonetheless deep.
Since basic tropical arithmetic is so simple, it occurs everywhere and thus has been rediscovered independently many times, long before tropical mathematics was recognized as a subject of its own.
However, it is the link to geometry, in its many flavors, which reveals the strength of these ideas.
The article on \emph{tropical convexity} by Develin and Sturmfels \cite{DevelinSturmfels04} provides ample evidence to support these claims.

Tropical convexity may be seen as the linear algebra part of tropical geometry.
It has a considerable pre-tropical history under the name of \enquote{$(\max,+)$-linear algebra}; see, e.g., the monographs of Cuninghame-Green \cite{Cuninghame-Green79}, Baccelli, Cohen, Olsder and Quadrat \cite{BaccelliCohenOlsderQuadrat92} and Butkovi\v{c} \cite{Butkovic:2010}.
The DNA of $(\max,+)$-linear algebra comprises quite a diverse range of topics including combinatorial optimization, functional analysis, discrete event systems and game theory.
One way of seeing tropical geometry is as a method to employ techniques from polyhedral geometry to solve problems in algebra and algebraic geometry.
As its fundamental contribution the article \cite{DevelinSturmfels04} fuses $(\max,+)$-linear algebra with that train of ideas.
The benefit is mutual:
Principles from algebraic geometry become accessible to the above named fields of application.
Conversely, these applications help to identify interesting geometric questions.

The structure of this survey is as follows.
In Section~\ref{sec:tconv} we start out with first principles concerning \emph{tropical cones} and \emph{tropical polytopes}, and we finish with the main result of \cite{DevelinSturmfels04}.
Throughout we occasionally deviate from the original notation and terminology of \cite{DevelinSturmfels04}.
While some of the original results allow for natural generalizations (which were found later), here we choose a path which is meant to ease the reading.
Section~\ref{sec:tlinear} is devoted to \emph{tropical linear spaces}, which can be defined and motivated in many ways; here we describe them as a special class of tropical polytopes.
The final Section~\ref{sec:applications} comprises examples showing the impact of tropical convexity outside tropical geometry.

For more details on tropical geometry in general we refer the reader to the monograph by Maclagan and Sturmfels \cite{Tropical+Book} and to my own book \cite{ETC}.
I am indebted to George Balla and two anonymous reviewers who helped me in improving the exposition.

\section{Tropically Convex Sets and Their Structure}
\label{sec:tconv}
For tropical arithmetic, as in \cite{DevelinSturmfels04}, we use the notation $\alpha \oplus \beta=\min(\alpha,\beta)$ and $\alpha\odot\beta=\alpha+\beta$ if $\alpha,\beta\in\TT:=\RR\cup\{\infty\}$.
The triplet $(\TT,\oplus,\odot)$ is the \emph{tropical semiring} (with respect to $\min$).
Tropical addition and scalar multiplication can be extended to vectors componentwise by letting $x\oplus y=(x_i\oplus y_i)_{i=1}^d$ and $\lambda\odot x=(\lambda \odot x_i)_{i=1}^d$, where $x,y\in\TT^d$ and $\lambda\in\TT$.

\subsection{Tropical cones}
The definition which marks the beginning of our topic is the following.
A set $S\subset\RR^d$ is called a \emph{tropical cone} if
\begin{equation}\label{eq:tropical-cone}
  \lambda\odot x \,\oplus\, \mu \odot y\;\in\; S \quad \text{for all } x,y\in S \text{ and } \lambda,\mu\in\RR \enspace .
\end{equation}
Tropical addition is idempotent, and so it immediately follows that if a tropical cone $S$ contains a vector $x\in\RR^d$, then it also contains $\smallSetOf{\lambda\odot x}{\lambda\in\RR}=x+\RR\1$, where $\1$ is the all-ones-vector (here of length $d$).
The following is \cite[Proposition 4]{DevelinSturmfels04}; 
see also \cite[Lemma 5.7]{ETC}.
\begin{proposition}
  A subset $S$ of $\RR^d$ is a tropical cone if and only if
  \[
    S \ = \ \SetOf{ \lambda_1 \odot x^{(1)} \,\oplus\, \dots \,\oplus\, \lambda_k \odot x^{(k)}}{\lambda_i\in\RR,\, x^{(i)}\in S } \enspace .
  \]
\end{proposition}
Tropical cones are precisely the \enquote{idempotent semimodules}, with respect to the tropical semiring, in the terminology of Cohen, Gaubert and Quadrat \cite{CohenGaubertQuadrat05}; see also \cite{Gaubert:PhD}.

\begin{figure}[bt]
  \newcommand{\scaling}{0.85}
  \centering
  \footnotesize
  \begin{tikzpicture}[scale=\scaling]
    \StandardGitter
    \draw[Poly] (0,0)--(-1,-1);
    \draw[Poly] (0,0)--(1,0);
    \draw[Poly] (0,0)--(0,1);
    \node at (-1.5,1.5) {$(a)$};
  \end{tikzpicture}
  \hfill
  \begin{tikzpicture}[scale=\scaling]
    \StandardGitter
    \draw[Poly] (-1,-1) -- (0,-1) -- (1,0) -- (1,1) -- (0,1) -- (-1,0) -- cycle;
    \node at (-1.5,1.5) {$(b)$};
  \end{tikzpicture}
  \hfill
  \begin{tikzpicture}[scale=\scaling]
    \StandardGitter
    \draw[Poly] (-1.2,-1.2) -- (1.2,-1.2) -- (1.2,1.2) -- cycle;
    \node at (-1.5,1.5) {$(c)$};
  \end{tikzpicture}

  \medskip
  
  \begin{tikzpicture}[scale=\scaling]
    \StandardGitter
    \draw[Poly,draw=none] (-1.7,-1.7) -- (0,0) -- (0,1.7) -- (1.7,1.7) -- (1.7,-1.7) -- cycle;
    \draw[Poly,fill=none] (-1.7,-1.7) -- (0,0) -- (0,1.7);
    \node at (-1.5,1.5) {$(d)$};
  \end{tikzpicture}
  \hfill
  \begin{tikzpicture}[scale=\scaling]
    \StandardGitter
    \draw[Poly] (-1.2,-1.2) -- (1.2,-1.2) -- (-1.2,0.8) -- cycle;
    \node at (-1.5,1.5) {$(e)$};
  \end{tikzpicture}
  \hfill
  \begin{tikzpicture}[scale=\scaling]
    \StandardGitter
    \draw[Poly] (-0.6,-0.6) -- (1.2,-0.6) arc [start angle=0, end angle=90, radius=1.8] -- cycle;
    \draw[Poly] (-1.3,-1.3) -- (-0.6,-0.6);
    \node at (-1.5,1.5) {$(f)$};
  \end{tikzpicture}

  \caption{Six tropically convex sets in the plane $\torus{3}$.
    Here and below a point $x\in\RR^3$ is drawn in the Euclidean plane $\RR^2$ as $(x_2-x_1,x_3-x_1)$; i.e., we consistently pick $(0,x_2-x_1,x_3-x_1)$ as a representative of $x+\RR\1$.
  }
  \label{fig:six-sets}
\end{figure}

A subset of the \emph{tropical projective torus} $\torus{d}$ is \emph{tropically convex} if it arises as the canonical projection $S/\RR\1$ of a tropical cone $S\subset\RR^d$.
For $M\subset\torus{d}$ the \emph{tropical convex hull}, written $\tconv(M)$, is the smallest tropically convex set containing $M$.
The most influential concept from \cite{DevelinSturmfels04} until today is the following:
A \emph{tropical polytope} is the tropical convex hull of finitely many points.
Figure \ref{fig:six-sets} shows six tropically convex sets in the plane, and the three examples in the top rows are finitely generated, i.e., tropical polytopes.
Yet the bottom three are not finitely generated.
For instance, the unique minimal set of generators of the triangle (e) is the uncountable set of points on the diagonal edge.
A first observation reveals a similarity between tropical and ordinary polytopes; see \cite[Proposition 21]{DevelinSturmfels04}.
The statement can be generalized to arbitrary compact tropically convex sets \cite[Theorem 5.43]{ETC}.
\begin{proposition}
  For each tropical polytope $P$ in $\torus{d}$ there exists a unique minimal set $V$ such that $P = \tconv(V)$.
\end{proposition}
The points in this minimal generating set are the \emph{tropical vertices} of $P$.
It is interesting that tropical cones and tropical polytopes occurred in the literature many years before, e.g., in Gaubert's thesis \cite{Gaubert:PhD}.
Yet it is the decidedly combinatorial point of view of \cite{DevelinSturmfels04} which reveals crucial new structural information.

Before we explore these structural results, it is worth spending a moment to confirm that the notions \enquote{(tropical) cone} and \enquote{(tropical) convexity} are adequate.
This is best seen by switching to the $\max$-tropical semiring $(\RR\cup\{-\infty\},\max,+)$, which is isomorphic to the $\min$-tropical semiring via passing to negatives.
Then any real number $\alpha\in\RR$ satisfies $\alpha>-\infty$, and $-\infty$ is the neutral element of the $\max$-tropical addition.
In this sense a tropical cone is closed with respect to taking tropical linear combinations with coefficients which are automatically positive.
This \enquote{built-in positivity} is responsible for the hybrid character of tropically convex sets: as we will see below, to some extent they behave like tropical analogs of linear spaces, and to some extent they behave like tropical analogs of convex sets.
Loho and V\'egh developed a generalization of tropical convexity which takes arbitrary (suitably defined) tropical signs into account \cite{LohoVegh:1906.06686}.

\begin{figure}[ht]

  \begin{tikzpicture}[scale=1]

    \tikzstyle{generator} = [circle, scale=0.6, draw=strawberry, fill=white]
    \tikzstyle{vertex} = [rectangle, scale=0.8, draw=strawberry, fill=white]
    
    \begin{scope}
      \clip (-3,-2) rectangle (6,3);

      \foreach \x/\y in {1/0, 1/2, -2/-1, 3/0, 4/2, 5/1}{
        \DrawMaxTropLine{Poly}{\x}{\y}{9};
      }
        
      \foreach \x/\y in {1/0, 1/2}{
        \node at (\x,\y) [generator] {};
      }
      
      \foreach \x/\y in {-2/-1, 3/0, 4/2, 5/1}{
        \node at (\x,\y) [generator] {};
      }
    \end{scope}

    \PolygonGitter
    
  \end{tikzpicture}
  
  \caption{Max-tropical hyperplane arrangement in $\torus{3}$.}
  \label{fig:arrangement}
\end{figure}

In the sequel we want to use the $\min$- and $\max$-tropical semirings together.
Whenever we need to distinguish these two tropical additions, we write $\minoplus=\min$ and $\maxoplus=\max$.
Roughly speaking, the $\min$-convention sides well with standard conventions concerning Puiseux series, while tropicalization with respect to $\max$ preserves the order relation; so the latter is the usual choice in optimization.
Conceptually it is clearly possible to express all tropical geometry in terms of either $\min$ or $\max$.
However, as in \cite{ETC}, we prefer to explain $\min$-tropical convexity in terms of $\max$-tropical hyperplanes.
This is the topic of the subsequent section.

\subsection{Arrangements of tropical hyperplanes}
For $d=1$ the tropical projective torus $\torus{1}$ is a single point, which is not very interesting.
So we assume that $d\geq 2$.
In the sequel we will fix finitely many vectors $v^{(1)},\dots,v^{(n)}\in \RR^d$, which we may identify with the columns of a matrix $V\in\RR^{d\times n}$.
We want to study the tropical polytope $P=\tconv(v^{(1)},\dots,v^{(n)})$, which we also write as $\tconv(V)$.
To this end it turns out to be instrumental that the matrix $V=(v_{ij})$ can be employed in more than one way.
For instance, the vector $-v^{(j)}$ gives rise to the homogeneous max-tropical linear form
\begin{equation}\label{eq:max-linear-form}
  -v_{1j} \odot X_1 \,\maxoplus\, \cdots \,\maxoplus\, -v_{dj} \odot X_d
\end{equation}
in the $d$ unknowns $X_1,\dots,X_d$.
The \emph{($\max$-)tropical hyperplane} $\cT(-v^{(j)})$ is the locus in $\torus{d}$, where this tropical polynomial vanishes, i.e., the maximum in \eqref{eq:max-linear-form} is attained at least twice.
It is immediate that the point $v^{(j)}$ is contained in $\cT(-v^{(j)})$, since in that case the evaluation at each term gives zero; as $d\geq 2$ that zero (which is also the maximum) is attained at least twice.
The point $v^{(j)}$ is the \emph{apex} of $\cT(-v^{(j)})$.
Recall that we defined tropical convexity in terms of $\min$ as tropical addition, whereas \eqref{eq:max-linear-form} is a $\max$-tropical polynomial.

\begin{figure}[ht]

\begin{tikzpicture}[x  = {(1cm,0cm)},
                    y  = {(0cm,1cm)},
                    z  = {(0cm,0cm)},
                    scale = 1,
                    color = {lightgray}]

  \PolygonGitter
                    
  \coordinate (v0_1) at (-2, -1);
  \coordinate (v1_1) at (-1, 0);

  \draw[Poly] (v1_1) -- (v0_1);

  \coordinate (v0_2) at (4, 1);
  \coordinate (v1_2) at (5, 1);

  \draw[Poly] (v1_2) -- (v0_2);

  \coordinate (v0_3) at (1, 0);
  \coordinate (v1_3) at (0, 1);
  \coordinate (v2_3) at (1, 1);
  \coordinate (v3_3) at (-1, 0);

  \draw[Poly] (v3_3) -- (v0_3) -- (v2_3) -- (v1_3) -- (v3_3) -- cycle;

  \coordinate (v0_4) at (1, 0);
  \coordinate (v1_4) at (1, 1);
  \coordinate (v2_4) at (2, 1);

  \draw[Poly] (v0_4) -- (v2_4) -- (v1_4) -- (v0_4) -- cycle;

  \coordinate (v0_5) at (1, 0);
  \coordinate (v1_5) at (2, 1);
  \coordinate (v2_5) at (3, 0);
  \coordinate (v3_5) at (4, 1);

  \draw[Poly] (v0_5) -- (v2_5) -- (v3_5) -- (v1_5) -- (v0_5) -- cycle;

  \coordinate (v0_6) at (0, 1);
  \coordinate (v1_6) at (1, 1);
  \coordinate (v2_6) at (1, 2);

  \draw[Poly] (v1_6) -- (v2_6) -- (v0_6) -- (v1_6) -- cycle;

  \coordinate (v0_7) at (1, 1);
  \coordinate (v1_7) at (2, 1);
  \coordinate (v2_7) at (1, 2);
  \coordinate (v3_7) at (3, 2);

  \draw[Poly] (v0_7) -- (v1_7) -- (v3_7) -- (v2_7) -- (v0_7) -- cycle;

  \coordinate (v0_8) at (2, 1);
  \coordinate (v1_8) at (4, 2);
  \coordinate (v2_8) at (4, 1);
  \coordinate (v3_8) at (3, 2);

  \draw[Poly] (v0_8) -- (v2_8) -- (v1_8) -- (v3_8) -- (v0_8) -- cycle;

  \tikzstyle{pseudovertex} = [circle, scale=0.6, fill=strawberry]
  \tikzstyle{generator} = [circle, scale=0.6, draw=strawberry, fill=white]
  \tikzstyle{vertex} = [rectangle, scale=0.8, draw=strawberry, fill=white]
  
  \foreach \x/\y in {0/1, 1/1, -1/0, 2/1, 4/1, 3/2}{
    \node at (\x,\y) [pseudovertex] {};
  }

  \foreach \x/\y in {1/0, 1/2}{
    \node at (\x,\y) [generator] {};
  }

  \foreach \x/\y in {-2/-1, 3/0, 4/2, 5/1}{
    \node at (\x,\y) [vertex] {};
  }

\end{tikzpicture}

  \caption{Min-tropical polygon in $\torus{3}$.
    The six generators are white; among these the four tropical vertices are marked by squares.
    The two circled generators are redundant for the tropical convex hull.
  }
  \label{fig:polygon}
\end{figure}
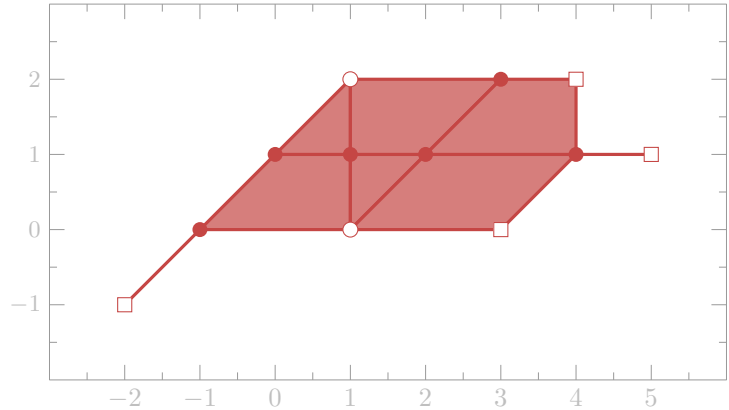

Let $A=A_V$ be the product of the $n$ linear forms $-v^{(1)}, \dots,-v^{(n)}$ in the semiring of homogeneous $\max$-tropical polynomials.
Then the $\max$-tropical hypersurface
\[
  \cT(A) \ = \ \cT(-v^{(1)}) \cup \cdots \cup \cT(-v^{(n)})
\]
is the union of $n$ tropical hyperplanes in $\torus{d}$.
The topological closure of a connected component of the complement $\torus{d}\setminus\cT(A)$ is called a \emph{region} of $A$.
Each region is an ordinary convex polyhedron in $\torus{d}\simeq\RR^{d-1}$, and the regions form a polyhedral subdivision of the tropical projective torus; we denote the latter as $\cS_V$.

A version of the \emph{Structure Theorem of Tropical Convexity} is the following.
Its essence is already expressed in \cite[Theorems 1 and 15, Proposition 16]{DevelinSturmfels04}.
However, the connection to tropical hyperplane arrangements was made explicit by Ardila and Develin \cite{ArdilaDevelin07}.
\begin{theorem}\label{thm:structure}
  For $V\in\RR^{d\times n}$ we have:
  \begin{enumerate}
  \item The polyhedral decomposition $\cS_V$ of $\torus{d}$, which is formed by the regions of the $\max$-tropical hyperplane arrangement $A_V$, is dual to the (lower) regular subdivision induced by $V$, considered as a height function on the vertices of the ordinary polytope $\Delta_{d-1}\times\Delta_{n-1}$.
  \item\label{it:tropical-complex} The $\min$-tropical polytope $\tconv(V)$ agrees with the union of the bounded cells of the polyhedral complex $\cS_V$.
  \end{enumerate}
\end{theorem}

In \cite{DevelinSturmfels04} the decomposition of $\tconv(V)$ as an ordinary polyhedral complex from Theorem~\ref{thm:structure}\eqref{it:tropical-complex} is called the \enquote{tropical complex} induced by $V$.
A more modern term for the same is \emph{covector decomposition}; this refers to Ardila and Develin \cite{ArdilaDevelin07}, who called a \enquote{covector} what was called a \enquote{type} in  \cite{DevelinSturmfels04}.
The covector decomposition is linearly isomorphic with the subcomplex of bounded faces of the unbounded ordinary polyhedron dual to $\cS_V$.

\begin{example}\label{exmp:polygon}
  We consider the $3{\times}6$-matrix
  \[
    V \ = \ \begin{pmatrix}
0 & 0 & 0 & 0 & 0 & 0\\
-2 & 1 & 3 & 4 & 5 & 1\\
-1 & 0 & 0 & 2 & 1 & 2
    \end{pmatrix} \enspace ,
  \]
  whose columns form the apices of six max-tropical lines in $\torus{3}$.
  We obtain an arrangement $A_V$ with 20 regions; see Figure~\ref{fig:arrangement}.
  Six of these regions are bounded, and together with two bounded edges, they form the min-tropical polytope $\tconv(V)$; see Figure~\ref{fig:polygon}.
\end{example}

Various extensions of Theorem~\ref{thm:structure} are known.
Notably, Ardila and Develin \cite{ArdilaDevelin07} introduced \enquote{tropical oriented matroids} to generalize tropical hyperplane arrangements.
Horn \cite{Horn:FPSAC:2012} showed that tropical oriented matroids correspond bijectively to the subdivisions of products of simplices which are not necessarily regular.
In Section \ref{subsec:infinite} we will see how infinite coordinates lead to \emph{subpolytopes} of $\Delta_{d-1}\times\Delta_{n-1}$ (i.e., convex hulls of subsets of the vertices) and their subdivisions.



\subsection{Products of simplices and their triangulations}

Let $V\in\RR^{d\times n}$.
The Structure Theorem \ref{thm:structure} ties the tropical convex hull $\tconv(V)$ to the regular subdivision $\cS_V$ of $\Delta_{d-1}\times\Delta_{n-1}$.
If $\cS_V$ is a triangulation, the matrix $V$ and the tropical polytopes $\tconv(V)$ are called \emph{generic}.
The product $S_d\times S_n$ of symmetric groups operates naturally on the vertices of $\Delta_{d-1}\times\Delta_{n-1}$.
For generic $V,W\in\RR^{d\times n}$ the two triangulations $\cS_V$ and $\cS_W$, are \emph{equivalent} if they lie in the same orbit with respect to the induced action.
It is an interesting task, but difficult, to enumerate these equivalence classes, also called \emph{combinatorial types}, for given $d$ and $n$.

\begin{example}\label{exmp:prism}
  For $d=2$ the product $\Delta_{d-1}\times\Delta_{n-1}$ is a prism over $\Delta_{n-1}$.
  This polytope has exactly $n!$ labeled triangulations, which form a single orbit; see \cite[\S6.2.1]{Triangulations} and \cite[\S4.6]{ETC}.
  So the first row of Table~\ref{tab:products} below contains ones.
\end{example}

Standard software for enumerating (regular) triangulations up to symmetry includes \TOPCOM \cite{TOPCOM-paper,TOPCOM} and \mptopcom \cite{JordanJoswigKastner:2018,mptopcom}.
The known values for products of simplices are given in Table~\ref{tab:products}.
The corresponding table in \cite[p.20]{DevelinSturmfels04} contains two typos, which have been corrected in \cite[p.241]{Tropical+Book}.
The first nontrivial values were probably computed by a \Maple program called \texttt{PUNTOS}, which was written by Jes\'us De Loera.
The table \cite[p.241]{Tropical+Book} was obtained with \TOPCOM, and the additional values, for $(d,n)\in\{(3,7), (4,5)\}$, have been found with \mptopcom \cite{JordanJoswigKastner:2018}.
With \TOPCOM v1.1.0 Rambau recently computed the total number of symmetry classes of all triangulations, including the nonregular ones, of $\Delta_3\times\Delta_5$ as 25{,}606{,}173{,}722; yet the corresponding regular count seems to be currently unknown.

\begin{table}[th]\centering
  \caption{The known numbers of combinatorial types of regular triangulations of $\Delta_{d-1}\times\Delta_{n-1}$, up to symmetry.}
\label{tab:products}
\begin{tabular}{crrrrrrr}
  \toprule
  $d$ $\backslash$ $n$ & 2 & 3 & 4 & 5 & 6 & 7 \\
  \midrule
  2 & 1 & 1 & 1 & 1 & 1 & 1 \\
  3 & & 5 & 35 & 530 & 13\,621 & 531\,862 \\
  4 & & & 7\,869 & 7\,051\,957 \\
  \bottomrule
\end{tabular}
\end{table}

\begin{remark}
  The product of two simplices $\Delta_{d-1}\times\Delta_{n-1}$ is known to admit nonregular triangulations if and only if $(d - 2)(n - 2) \geq 4$; see \cite[Theorem 6.2.19]{Triangulations}.
  While such nonregular triangulations are not directly relevant for tropical geometry, their mere existence (and abundance, in higher dimensions) is one of the obstacles to overcome in obtaining full enumerations of the regular triangulations.
\end{remark}

\subsection{Infinite coordinates}
\label{subsec:infinite}

In the matrices $V$ that we considered so far the neutral element, $\infty$, of the tropical addition $\minoplus=\min$ was explicitly forbidden.
Let the \emph{support} of $V\in\TT^{d\times n}$ be the subset of $[d]\times[n]$ chosen by the finite coefficients.
Fink and Rinc\'on generalized part~(1) of the Structure Theorem \ref{thm:structure} to arbitrary $V$, where $\Delta_{d-1}\times\Delta_{n-1}$ is replaced by the subpolytope defined by the support of $V$; see \cite[\S4]{FinkRincon:2015} and \cite[Corollary 4.13]{ETC}.
In order to understand how part~(2) can be generalized, we first need to re-examine our basic definitions.

Clearly, we can replace $\RR^d$ by $\TT^d$ to obtain tropical cones in $\TT^d$, with the analogous definition to \eqref{eq:tropical-cone}.
Let
\[
  \TP^{d-1} \ := \ \bigl( \TT^d\setminus\{\infty\1\} \bigr) \,/\, \RR\1
\]
be the \emph{tropical projective space} of rank $d-1$, where $\infty\1$ is the vector of length $d$ where each coefficient is $\infty$.
Now passing from $\TT^d\setminus\{\infty\1\}$ to $\TP^{d-1}$ by taking the quotient provides notions of tropical convexity and tropical polytopes in tropical projective spaces.

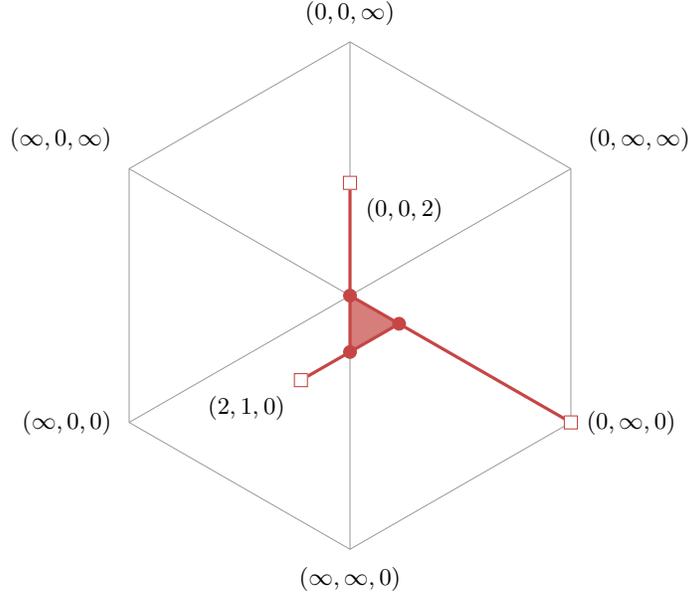
\begin{figure}[bt]
  \newcommand\localinf{4.5}
  
  \begin{tikzpicture}[x  = {(-0.87cm,-0.5cm)},
                    y  = {(0.87cm,-0.5cm)},
                    z  = {(0cm,1cm)},
                    scale = 0.75,
                    color = {lightgray}]


  \coordinate (e1) at (-\localinf,0,0);
  \coordinate (e2) at (0,-\localinf,0);
  \coordinate (e3) at (0,0,-\localinf);
  
  \coordinate (f1) at (\localinf,0,0);
  \coordinate (f2) at (0,\localinf,0);
  \coordinate (f3) at (0,0,\localinf);

  \draw[GitterStyle] (e1) -- (f3) -- (e2) -- (f1) -- (e3) -- (f2) -- cycle;

  \foreach \i in {1,2,3}{
    \draw[GitterStyle] (e\i) -- (f\i);
  }

  \coordinate (a) at (0,\localinf,0);
  \coordinate (b) at (0,0,2);
  \coordinate (c) at (0,-1,-2);

  \coordinate (u) at (0,1,0);
  \coordinate (v) at (0,0,0);
  \coordinate (w) at (0,0,-1);
  
  \filldraw[Poly] (u) -- (v) -- (w) -- cycle;
  \draw[Poly] (a) -- (u);
  \draw[Poly] (b) -- (v);
  \draw[Poly] (c) -- (w);

  \small 
  \node at (a) [vertex, label={[black] right:$(0,\infty,0)$}] {};
  \node at (b) [vertex, label={[black] below right:$(0,0,2)$}] {};
  \node at (c) [vertex, label={[black] below left:$(2,1,0)$}] {};

  \node at (u) [pseudovertex] {};
  \node at (v) [pseudovertex] {};
  \node at (w) [pseudovertex] {};

  \node at (f1) [label={[black] left:$(\infty,0,0)$}] {};
  \node at (f3) [label={[black] above:$(0,0,\infty)$}] {};

  \node at (e1) [label={[black] above right:$(0,\infty,\infty)$}] {};
  \node at (e2) [label={[black] above left:$(\infty,0,\infty)$}] {};
  \node at (e3) [label={[black] below:$(\infty,\infty,0)$}] {};
  
  \end{tikzpicture}
  \caption{Tropical triangle in $\TP^2$ (isometric sketch).  Here $\TP^2$ takes the shape of a regular hexagon, and $\torus{3}$ is its interior.}
  \label{fig:projective-triangle}
\end{figure}

The construction above is beneficial from a topological point of view as $\TP^{d-1}$ compactifies $\torus{d}$.
In fact, the pair $(\TP^{d-1},\torus{d})$ is homeomorphic with the pair formed by a $(d{-}1)$-dimensional closed ball and its interior, see \cite[Proposition 5.3]{ETC}; the planar case $d=3$ is sketched in Figure~\ref{fig:projective-triangle}.
The boundary equals the set difference $\partial\TP^{d-1}=\TP^{d-1} \setminus (\torus{d})$, and this is naturally stratified by the supports of vectors in $\TT^d$.
Here the \emph{support} of $v\in\TT^d$ is the set of indices $i$ such that $v_i$ is finite.
For $Z\subsetneq[d]$ let $\TP(Z)$ be the \emph{stratum} of points whose support equals the complement $[d]\setminus Z$.
This set is a copy of $\torus{d-k}$, where $k$ is the cardinality of $Z$.
Since the strata in $\TP^{d-1}$ correspond bijectively to the proper subsets of $[d]$, the tropical projective space $\TP^{d-1}$ inherits the combinatorial structure of a simplex.

\begin{observation}
  Let $P \subseteq \TP^{d-1}$ be a tropical polytope.
  For any $Z\subsetneq[d]$ the intersection $P\cap\TP(Z)$ is tropically convex.
\end{observation}

We call $P\cap\TP(Z)$ a \emph{stratum} of $P$.
In this way, by applying Theorem~\ref{thm:structure}(1) to the strata and taking suitable unbounded hyperplane regions into account, we arrive at a generalization of Theorem~\ref{thm:structure}(2) to tropical polytopes in tropical projective space. 
The following example shows that the strata of a tropical polytope are not necessarily finitely generated.

\begin{example}\label{exmp:projective-triangle}
  We consider the matrix
  \[
    V \ = \ \begin{pmatrix}
      0 & 0 & 2\\
      \infty & 0 & 1 \\
      0 & 2 & 0
    \end{pmatrix} \ \in \TT^{3\times 3} \enspace ,
  \]
  and $T=\tconv(V)$ is the tropical triangle in $\TP^2$ shown in Figure~\ref{fig:projective-triangle}.
  For instance, the stratum $\TP(\{2\})$ is the set of points with homogeneous coordinates $(v_1,\infty,v_3)$ where $v_1,v_3\in\RR$; the intersection $P\cap\TP(\{2\})$ comprises the single point $(0,\infty,0)$.
  All other boundary strata of $P$ are empty.
  The generic stratum $P\cap\TP(\emptyset)$ is tropically convex, but not a tropical polytope in the tropical projective torus $\torus{3}=\TP(\emptyset)$.
\end{example}

\subsection{Eigenvalues and eigenvectors}
\label{sec:eigen}
We start out with studying another tropical analog of ordinary linear algebra.
Let $A=(a_{ij})$ be an $n{\times}n$-matrix with coefficients in the tropical semiring $\TT$.
A real number $\lambda$ is a \emph{tropical eigenvalue} of $A$ if there is a vector $x\in\TT^n\setminus\{\infty\1\}$ such that
\begin{equation}
  \label{eq:eigen}
  A \odot x \ = \lambda \odot x \enspace ,
\end{equation}
in which case $x$ is a \emph{tropical eigenvector} of $A$ with respect to the eigenvalue $\lambda$.
One can read the matrix $A$ as the weighted adjacency matrix of a directed graph $\Gamma=\Gamma(A)$ on $n$ nodes.
Here $a_{ij}=\infty$ signals the absence of the arc $(i,j)$.
However, when $a_{ij}$ is finite, this value is the weight of that arc.
The \emph{mean weight} of a directed path $\pi$ with $k$ arcs $(i_0,i_1),(i_1,i_2),\dots,(i_{k-1},i_k)$ is the value $c(\pi)=\tfrac{1}{k}(a_{i_0,i_1}+\cdots+a_{i_{k-1},i_k})$.
For $i_k=i_0$ the path $\pi$ is a cycle, and the mean weight $c(\pi)$ is called the \emph{cycle mean} of $\pi$.
Now the \emph{minimum cycle mean} of $A$ is defined as
\[
  \lambda(A) \ = \min\SetOf{ c(\zeta) }{ \zeta \text{ cycle in } \Gamma(A) } \enspace .
\]
\begin{theorem}
  Let $\Gamma(A)$ be strongly connected.
  Then the minimum cycle mean $\lambda(A)$ is the only tropical eigenvalue of $A$.
\end{theorem}
Karp~\cite{Karp78} gave a polynomial time algorithm for computing $\lambda(A)$; see also \cite[\S3.6]{ETC}.
This tropical eigenvalue is also the smallest root of its characteristic tropical polynomial \cite[Corollary 5.1.5]{Tropical+Book}.

For the rest of this section we assume that $\Gamma(A)$ is strongly connected and $\lambda(A)$ equals zero.
In this case it is easy to describe the tropical eigenvectors of $A$.
Namely, we can form the \emph{Kleene star}, which is the matrix
\begin{equation}\label{eq:kleene}
  A^* \ = \ I \oplus A \oplus A^{\odot 2} \oplus \cdots \oplus A^{\odot(n-1)} \enspace,
\end{equation}
where $A^{\odot k}$ is the $k$th tropical matrix power, and $I=A^{\odot 0}$ is the tropical identity matrix, with zero diagonal and all off-diagonal entries equal to $\infty$.
Then we have $A\odot A^*=A^*$.
That is to say that each column of $A^*$ is a tropical eigenvector of $A$ with respect to the eigenvalue zero; see \cite[Proposition 3.39]{ETC}.
More generally, the \emph{tropical eigenspace} is the set
\[
  \Eig(A) \ = \ \SetOf{ x\in\TT^n\setminus\{\infty\1\} }{ A \odot x = x } \enspace .
\]
For $x,y\in\Eig(A)$ and $\lambda,\mu\in\RR$ we compute
\[
  A\odot(\lambda\odot x \oplus \mu\odot y) \ = \ \lambda\odot(A\odot x)\oplus \mu\odot(A\odot y) \ = \ \lambda\odot x\oplus \mu\odot y \enspace ,
\]
and this says that the tropical eigenspace $\Eig(A)$ is a tropical cone.
In fact, we have $\Eig(A)=\tconv(A^*)$ and, by construction, $\tconv(A^*)\subseteq\tconv(A)$.
\begin{example}
  The matrix $V$ from Example~\ref{exmp:projective-triangle} has minimum cycle mean $\lambda(V)=0$, and the Kleene star reads 
  \[
    V^* \ = \ \begin{pmatrix}
      0 & 1 & 0\\
      0 & 0 & 0\\
      1 & 1 & 0
    \end{pmatrix} \enspace .
  \]
  The directed graph $\Gamma(V)$ is strongly connected, and each column of $V^*$ lies in a cycle of weight zero.
  The tropical eigenspace is the tropical cone in $\TP^ {n-1}$ spanned by the columns of $V^*$.
  In Figure~\ref{fig:projective-triangle} the eigenspace $\Eig(V)=\Eig(V^*)$, or rather its image in $\TP^2$, is the unique triangular cell in the covector decomposition of $\tconv(V)$.
\end{example}

The tropical convex hulls of Kleene stars play a very special role in tropical convexity.
We continue with our general assumptions that $A\in\TT^{n\times n}$ is a matrix such that $\Gamma(A)$ is strongly connected.
Then $A^*$ has finite coefficients only, and  $\tconv(A^*)/\RR\1$ is a tropical polytope in the tropical projective torus.
We can identify $\torus{n}$ with $\RR^{n-1}$ via $(x_1,x_2,\dots,x_n)+\RR\1\mapsto (x_2-x_1,\dots,x_n-x_1)$; this identification was already used for drawing our figures.
The following statement is essentially equivalent to \cite[Proposition 18]{DevelinSturmfels04}; see also \cite[Theorem 6.38]{ETC}.
\begin{proposition}
  The tropical polytope $\tconv(A^*)$, seen as a subset of $\RR^{n-1}=\torus{n}$ is an ordinary convex polytope.
  We have \[ \tconv(A^*)=\SetOf{ x\in\torus{n} }{ x_i-x_j \leq a^*_{ij} } \enspace . \]
\end{proposition}
The tropical polytopes like $\tconv(A^*)$ are exactly those tropical polytopes which, seen as subsets of $\RR^{n-1}$, are convex in the ordinary sense.
In \cite{JoswigKulas:2010} they have been called \emph{polytropes}; see also \cite[\S6.5]{ETC}.
The bounded regions of tropical hyperplane arrangements, i.e., the bounded cells of the covector decomposition from the Structure Theorem~\ref{thm:structure} are exactly the polytropes.
Further, the integral polytropes, i.e., the tropical convex hulls of Kleene stars with integer coefficients, are precisely the \emph{alcoved polytopes} (of type A) of Lam and Postnikov \cite{LamPostnikov05}.
The hexagonal shape of the tropical projective plane in Figure~\ref{fig:projective-triangle} refers to the equality
\[
  \TP^2 \ = \ \tconv\begin{pmatrix} 0 & \infty & \infty \\ \infty & 0 & \infty \\ \infty & \infty & 0 \end{pmatrix} \enspace ,
\]
which means that we may view $\TP^{n-1}$ as an \enquote{infinitely large} tropical unit ball, which is a polytrope \cite[Example 6.39]{ETC}.

The Kleene star $A^*=(a^*_{ij})$ contains relevant information about the directed graph $\Gamma(A)$.
More precisely, the coefficient $a^*_{ij}$ records the length (or weight) of a shortest path from $i$ to $j$.
Computing shortest paths is a basic problem in combinatorial optimization.
Suitable algorithms are, e.g., at the core of common navigational devices.
In \cite{JoswigSchroeter:1904.01082} we studied a generalization, where the coefficients of the matrix $A$ are tropical polynomials.
Again the definition \eqref{eq:kleene} of the Kleene star makes sense, and then $A^*$ describes the lengths of shortest paths in a parametric setting.
The original motivation for \cite{JoswigSchroeter:1904.01082} was a new computational approach for classifying the combinatorial types of polytropes.
This continued previous work of Sturmfels--Tran \cite{SturmfelsTran:2013} and Tran \cite{Tran14}.
The parametric shortest path algorithms proposed in \cite{JoswigSchroeter:1904.01082} were used by Cleveland et al.~\cite{NASA} for developing delay tolerant networks, with the ultimate goal to construct the Solar System Internet (SSI).

\subsection{Tropical halfspaces}

As its main benefit the compactification into $\TP^{d-1}$ makes tropical convexity more symmetric in the following sense.
The main theorem for ordinary cones \cite[\S1.3 and \S1.4]{Ziegler:Lectures+on+polytopes} says that each ordinary cone is the intersection of finitely many linear halfspaces and conversely.
In the tropical projective space there is a corresponding result for tropical polytopes.
To describe it we need to explain what the tropical analogs to linear halfspaces are.

Let $a^+,a^-\in\TT^d\setminus\{\infty\1\}$.
For conciseness we assume that $\supp(a^+)\cap\supp(a^-)=\emptyset$.
Then the set $\smallSetOf{x\in\TT^d}{a^+\minodot x \leq a^-\minodot x}$ is a \emph{tropical halfspace} in $\TT^d$, and we use the same name for the image under the projection to $\TP^{d-1}$.
The following was proved by Gaubert and Katz \cite{GaubertKatz:2011}, and this corrects an error in \cite{Joswig:2005}.

\begin{theorem}\label{thm:double-description}
  Every tropical polytope in $\TP^{d-1}$ is the intersection of finitely many tropical halfspaces.
  Conversely, every such intersection is a tropical polytope.
\end{theorem}

While the above result suggests that tropical polytopes are similar to ordinary polytopes, there are remarkable distinctions.

\begin{example}
  The tropical triangle $T$ from Example~\ref{exmp:projective-triangle} is the intersection of five tropical halfspaces.
  More precisely, $T$ is the set of points $x\in\TP^{2}$ satisfying the homogeneous tropical matrix inequality
  \begin{equation}\label{eq:triangle-ineq}
    \begin{pmatrix}
      \infty & 0 & 0\\
      0 & \infty & 1\\
      0 & -1 & \infty\\
      \infty & \infty & -2\\
      -2 & \infty & \infty
    \end{pmatrix}
    \, \minodot \,
    x
    \ \leq \
    \begin{pmatrix}
      0 & \infty & \infty\\
      \infty & 0 & \infty\\
      \infty & \infty & 0\\
      0 & 0 & \infty\\
      \infty & -1 & 0
    \end{pmatrix} 
    \, \minodot \, 
    x
    \enspace .
  \end{equation}
  Each pair of corresponding rows of the two matrices to the left and to the right defines one tropical halfspace.
  The inequality description \eqref{eq:triangle-ineq} is irredundant.
  The tropical triangle $T$ is depicted in Figure~\ref{fig:projective-triangle}; yet the visualization in $\TP^2$ is not canonical in the sense that the boundary is deliberately drawn at some finite distance from the origin.
  Consequently, in such a picture some tropically convex sets would look nonlinear near the boundary; e.g., see Figure~\cite[Figure 6.5]{ETC}.
\end{example}

In ordinary polytope theory McMullen's Upper Bound Theorem \cite{McMullen:1970} dictates how many inequalities are necessary at most to describe a polytope; see also \cite[Theorem~8.23]{Ziegler:Lectures+on+polytopes}.
Allamigeon, Gaubert and Katz \cite{AllamigeonGaubertKatz:2011} proved an analog for tropical polytopes; see also \cite[Theorem 8.10]{ETC}.
The precise bounds differ slightly between the ordinary and the tropical setting, since nonnegativity constraints play a special role.
So it may happen that a tropical triangle, like the one in Example~\ref{exmp:projective-triangle}, requires five tropical linear inequalities.
Arrangements of tropical halfspaces can also be studied combinatorially.
This leads to yet another variation of the Structure Theorem \ref{thm:structure}; see \cite[\S7.4]{ETC}.

\subsection{Real Puiseux series}
\label{sec:real-puiseux}
After the discussion above it is natural ask: how exactly are ordinary and tropical polytopes related?
One way to answer this question goes through the tropicalization of ordinary polytopes defined over suitable ordered fields.
A \emph{(real) Puiseux series} is a formal power series of the form
\[
  \gamma(t) \ = \ \sum_{k=\ell}^\infty  c_k t^{k/N} \enspace ,
\]
where $\ell\in\ZZ$, $N\in\NN$, $c_k\in\RR$, and $c_\ell\neq 0$.
With the coefficient-wise addition and the usual product of formal power series the real Puiseux series form a field, which we denote by $\puiseuxseries{\RR}{t}$.
Moreover, this field is ordered: $\gamma(t)$ inherits the sign from its leading coefficient $c_\ell$.
So a real Puiseux series is nonnegative if it is zero or its leading coefficient is positive.
The orderd field $\puiseuxseries{\RR}{t}$ admits an additional structure, which becomes visible by considering the \emph{order map} $\val:\puiseuxseries{\RR}{t}\to\QQ\cup\{\infty\}$, which sends $\gamma(t)$ to its leading exponent $\ell/N$; and we set $\val(0)=\infty$.
The order map forms a discrete valuation.
We have
\[
  \val(\gamma + \delta) = \min\bigl(\val(\gamma),\val(\delta)\bigr) \quad \text{and} \quad \val(\gamma \cdot \delta) = \val(\gamma)+\val(\delta) \enspace ,
\]
provided that $\gamma$ and $\delta$ both are strictly positive.
This can be summarized as follows.

\begin{lemma}\label{lem:homomorphism}
  Restricting the order map to nonnegative real Puiseux series yields a homomorphism $\val:\puiseuxseries{\RR}{t}_{\geq 0} \to \TT$ of semirings.
\end{lemma}

A minor technical nuisance is that the order map is not surjective onto $\TT$.
This can be remedied in several ways, for instance by passing to a suitable extension field which is large enough.
One choice is a construction suggested by Markwig \cite{Markwig:2010}; see also \cite[\S2.7]{ETC}.
In the sequel we assume that $\LL$ is some real closed field \cite{BochnakCosteRoy:1998}, equipped with a valuation map, $\val$, which is surjective onto $\TT$.
For such fields Develin and Yu \cite{DevelinYu07} proved the following; see also \cite[Proposition 5.8]{ETC}.
Note that cones, polytopes, linear programming and such make sense over arbitrary ordered fields; e.g., see \cite{Jeroslow:1973b} and \cite[{\S}A.1]{ETC}.

\begin{theorem}\label{thm:develin-yu}
  For any ordinary polyhedral cone, $C$, which is contained in the nonnegative orthant $\LL_{\geq 0}^d$, the pointwise image $\val(C)$ is a tropical cone in $\TT^d$.
  Conversely, each tropical cone arises in this way.
\end{theorem}

This implies a similar result for ordinary and tropical polytopes by considering the quotients $(\LL_{\geq 0}^d\setminus\{0\})/\LL_{>0}$ and $\TP^{d}=( \TT^d\setminus\{\infty\1\} ) / \RR\1$.
The significance of Theorem~\ref{thm:develin-yu} comes from the fact that $\LL$ was assumed to be real closed.
In this case the Tarski--Seidenberg Transfer Principle \cite[Corollary 5.2.4]{BochnakCosteRoy:1998} entails that the ordered fields $\LL$ and $\RR$ are elementary equivalent.
In particular, the combinatorial types of polytopes over $\LL$ and over $\RR$ are the same; see \cite[Proposition 5.12]{Benchimol:2014}, \cite{JoswigLohoLorenzSchroeter:2016} and \cite[Theorem 8.36]{ETC}.

\begin{remark}\label{rem:order}
  For two real Puiseux series $\gamma,\delta\in\puiseuxseries{\RR}{t}$ with $\gamma \leq \delta$ we have $\val(\gamma)\geq\val(\delta)$, which reverses the ordering.
\end{remark}

\subsection{Further aspects of tropical convexity and beyond}

The connection between ordinary and tropical convexity can also be explored from a systematic point of view.
For instance, Cohen--Gaubert--Quadrat \cite{CohenGaubertQuadrat04,CohenGaubertQuadrat05} investigated subsemimodules and convex subsets of general idempotent semimodules.
To overcome the crucial nonnegativity assumption in Theorem \label{thm:develin-yu} Loho and V\'egh proposed a version of tropical convexity in the signed tropical semiring \cite{LohoVegh:2020}.
Motivated by aspects of game theory Akian--Gaubert--Vanucci introduced the concept of \enquote{ambitropical convexity} \cite{AkianGaubertVannucci:2108.07748}.
As a vast generalization in a different direction, Alessandrini~\cite{Alessandrini:2013} explored the tropicalization of arbitrary semialgebraic sets; Allamigeon--Gaubert--Skomra studied the special case of tropical spectrahedra \cite{AllamigeonGaubertSkomra:2020}.

\section{Tropical Linear Spaces}
\label{sec:tlinear}
Speyer and Sturmfels \cite{SpeyerSturmfels04} investigated the tropicalization of ordinary linear spaces and their parameter spaces.
This leads to tropical linear spaces, which arise as certain tropical polytopes in a tropical projective space.
We fix $d$ and $n$ such that $d\leq n$.
A \emph{finite tropical Plücker vector} is a map $\pi: \tbinom{[n]}{d} \to \RR$ such that the minimum of the three numbers
\[
  \pi(\rho+ij) + \pi(\rho+k\ell) \,,\   \pi(\rho+ik) + \pi(\rho+j\ell) \,,\   \pi(\rho+i\ell) + \pi(\rho+jk)
\]
is attained at least twice for any $(d{-}2)$-subset $\rho\subset[n]$ and pairwise distinct $i,j,k,\ell\in[n]-\rho$.
Here we use the shortcut $\rho+ij$ for the set $\rho\cup\{i,j\}$; more generally, we employ this notation for the disjoint union of sets where the second set is written as a list of its elements.
Let $\sigma\subset[n]$ have cardinality $d-1$.
Then $\sigma*$ is defined as the vector in $\TT^n$ whose $i$th coordinate is $\pi(\sigma+i)$ if $i\not\in\sigma$ and $\infty$ otherwise.
The vector $\sigma*$ is a \emph{cocircuit} of $\pi$.
Now the tropical convex hull in $\TP^{n-1}$ of all cocircuits of $\pi$ is the \emph{uniform tropical linear space} $\cL_\pi$.

\begin{figure}[tb]
  \begin{tikzpicture}[scale=2]
    \tikzstyle{pseudovertex} = [circle, scale=0.6, fill=strawberry]
    \tikzstyle{generator} = [circle, scale=0.6, draw=strawberry, fill=white]
    \tikzstyle{vertex} = [rectangle, scale=0.8, draw=strawberry, fill=white]

    \node[label=left:{\footnotesize $(0,7,1,2,3)$}] (v5) at (-1,-0.2) [pseudovertex] {};
    \node[label=below:{\footnotesize $(2,4,3,4,0)$}] (v6) at (0,0) [pseudovertex] {};
    \node[label=right:{\footnotesize $(3,4,4,4,0)$}] (v7) at (1,-0.2) [pseudovertex] {};
    \node (v8) at (2,0.5) [vertex] {1};
    \node (v9) at (2,-0.9) [vertex] {3};
    \node (v10) at (0,1) [vertex] {4};
    \node (v11) at (-2, -0.9) [vertex] {5};
    \node (v12) at (-2,0.5) [vertex] {2};

    \draw (v5) -- (v6) -- (v7);
    \draw (v10) -- (v6);
    \draw (v11) -- (v5) -- (v12);
    \draw (v8) -- (v7) -- (v9);
  \end{tikzpicture}
  \caption{A tropicalized linear space in $\TP^4$.
    Leaf labels correspond to row/column indices of the cocircuit matrix \eqref{eq:cocircuits}.}
  \label{fig:tropicalized}
\end{figure}
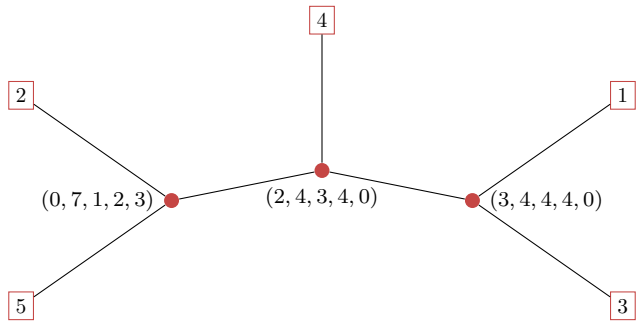

\begin{example}\label{exmp:tropicalized}
  With $d=2$ and $n=5$ we consider the $2{\times}5$-matrix
  \[
    M \ = \ \begin{pmatrix}
      1 & 0 & t & t^2 & t^3 \\
      0 & 1 & 1 & t-1 & t^{-4}
    \end{pmatrix} \enspace ,
  \]
  whose entries are Puiseux series.
  Denoting the $2{\times}2$ submatrix formed by columns $i$ and $j$ by $M^{(ij)}$, we get a map $\pi:\tbinom{[5]}{2}\to\RR$ by setting $\pi(ij)=\val(\det(M^{(ij)})$.
  This reads
  \[
    \begin{aligned}
      &{12} \mapsto 0 \,,\ {13} \mapsto 0 \,,\ {14} \mapsto 0 \,,\ {15} \mapsto -4 \,,\ {23} \mapsto 1 \,,\\
      &{24} \mapsto 2 \,,\ {25} \mapsto 3 \,,\ {34} \mapsto 1 \,,\ {35} \mapsto -3 \,,\ {45} \mapsto -2 \enspace .
    \end{aligned}
  \]
  Since the maximal minors of a matrix satisfy the Plücker relations, it follows that $\pi$ is a tropical Plücker vector; see \cite[\S4.3]{Tropical+Book}.
  The five cocircuits are the rows (or columns) of the symmetric matrix
  \begin{equation}\label{eq:cocircuits}
    \begin{pmatrix}
      \infty & 0 & 0 & 0 & -4\\
      0 & \infty & 1 & 2 & 3\\
      0 & 1 & \infty & 1 & -3\\
      0 & 2 & 1 & \infty & -2\\
      -4 & 3 & -3 & -2 & \infty
    \end{pmatrix} \enspace .
  \end{equation}
  Their tropical convex hull $\cL_\pi$ in $\TP^4$, intersected with $\torus{5}$, is the stratum $\cL\cap\TP(\emptyset)$; it is a contractible one-dimensional polyhedral complex with two bounded edges and five unbounded ones.
  They form a trivalent tree with five leaves, one for each cocircuit; the leaves are the generators, and they sit in the boundary $\partial\TP^4=\TP^4\setminus\torus{5}$.
  The stratum $\cL_\pi\cap\TP(\{i\})$ comprises the $i$th cocircuit only, and all other strata are empty.
  The tree is shown in Figure~\ref{fig:tropicalized}, with homogeneous coordinates in $\torus{5}$ for the three interior nodes.
\end{example}

A tropical linear space is called a \emph{tropicalized linear space} if it arises from the procedure followed in Example~\ref{exmp:tropicalized}.
In this case we start out with a matrix whose coefficients are complex Puiseux series, and the resulting tropical linear space is the tropicalization of the row space of that matrix.
Classically, the row spaces of $d{\times}n$-matrices of full row rank $d$ parameterize the $d$-dimensional linear subspaces in $n$-space, and they form the Grassmannian $\Gr(d,n)$.
The Grassmannian is an algebraic variety, defined over the field of complex Puiseux series, cut out by the Plücker relations.
However, since the coefficients of the Plücker relations are integers, i.e., complex Puiseux series which are constants, it suffices to view $\Gr(d,n)$ as a complex variety.
Tropicalizing $\Gr(d,n)$, we arrive at the tropical variety $\TGr(d,n)$, which is called the \emph{tropical Grassmannian}; see \cite{SpeyerSturmfels04} or \cite[\S4.3]{Tropical+Book}.
The tropical Grassmannian $\TGr(d,n)$ is a polyhedral fan in $\RR^{\tbinom{n}{d}}$.
That fan is pure, and all its maximal cones share the same dimension $(n-d)d+1$.
After modding out linealities and passing to the intersection with the unit sphere we arrive at a spherical polytopal complex of dimension $nd-n-d^2$.

\begin{theorem}[{Speyer and Sturmfels \cite{SpeyerSturmfels04}}]
  The points in $\TGr(d,n)$ are precisely the tropical Plücker vectors of tropicalized linear spaces.
\end{theorem}

We have $\TGr(d,n)\cong\TGr(n-d,n)$ as polyhedral fans.
Hence, for the rest of this section, we assume $n\geq 2d$.
For $d=2$ the tropical Grassmannians $\TGr(2,n)$ are spaces of phylogenetic trees on $n$ labeled leaves.
In this case each tropical linear space is tropicalized.
However, for $d\geq 3$ and $n\geq 7$ this is no longer true.
In these cases the \emph{Dressian} $\Dr(d,n)$, which is the polyhedral fan comprising all tropical Plücker vectors, is strictly larger than $\TGr(d,n)$.

\begin{remark}
  The tropical Grassmannians $\TGr(3,n)$ and the Dressians $\Dr(3,n)$ have explicitly been determined for the first few values of $n$: $n=6$ already occurs in \cite{SpeyerSturmfels04}, the case $n=7$ was resolved in \cite{HerrmannJensenJoswigSturmfels:2009}.
  The current state of the art is marked by $n=8$.
  Despite that the Dressian is (much) larger than the tropical Grassmannian it is easier to compute; this was accomplished already in \cite{HerrmannJoswigSpeyer14}.
  Computing $\TGr(3,8)$ is a very recent achievement, which required sophisticated technology and nontrivial parallel algorithms in computer algebra \cite{BendleBoehmRenSchroeter:2003.13752}.
  About $\TGr(4,8)$ and $\Dr(4,8)$ and their even bigger friends we only have very limited partial information.
\end{remark}

\section{Applications}
\label{sec:applications}
Recent years saw an increasing number of applications of tropical geometry, and tropical convexity especially, to many fields of mathematics and beyond.
Here we want to highlight one application to topics in computational complexity.
We close this section with a brief collection of further references.

\subsection{Linear programming}
An ordinary linear program is an optimization problem which seeks to minimize (or maximize) a linear function over an ordinary polyhedron.
Classically this has been studied over the reals, but the concept makes sense over any ordered field, in particular, the field of real Puiseux series studied in Section~\ref{sec:real-puiseux}.
Before we look at a tropical analog, we want to re-examine our setup.

By Remark~\ref{rem:order} the valuation map $\val$, which sends a real Puiseux series to its lowest degree, reverses the order.
If we want to relate ordinary and tropical inequalities, this change of orientation can be confusing.
To this end, we can instead take dual Puiseux series, which we denote by $\puiseuxseries{\RR}{t}^*$, and whose exponents form a descending sequence.
Then the \emph{dual valuation}, $\val^*$, takes a dual Puiseux series to its highest degree.
As in Lemma~\ref{lem:homomorphism} we obtain a semiring homomorphism $\val^*:\puiseuxseries{\RR}{t}^*_{\geq0}\to\TTmax$, where $\TTmax=(\RR\cup\{-\infty\},\max,+)$ is the max-tropical semiring.
However, $\val^*$ preserves the ordering.

So we consider two matrices $A^+,A^-\in\TTmax^{m\times n}$, two vectors $b^+,b^-\in\TTmax^m$ and another vector $c\in\TTmax^n$.
Often we abbreviate $A=(A^+,A^-)$ and $b=(b^+,b^-)$.
Equipped with this input, the \emph{tropical linear program} $\TLP(A,b,c)$ is the optimization problem
\begin{equation}\label{eq:tlp}
  \begin{array}{l@{\quad}l}
    \text{minimize}   & \tscp{c}{x} \\[\jot]
    \text{subject to} & A^+ \maxodot x \maxoplus b^+ \geq A^- \maxodot x \maxoplus b^-\\[\jot]
    & x\in\TTmax^{n} \enspace .
  \end{array}
\end{equation}
Each row of the constraints corresponds to one tropical halfspace as in Theorem~\ref{thm:double-description}, but there are two minor differences.
First, as pointed out above, here we consider tropical arithmetic with respect to $\oplus=\maxoplus=\max$ as the tropical addition.
This entails a tropical matrix multiplication which uses that version of tropical addition; to stress this we write the latter as \enquote{$\maxodot$}.
Second, by taking the vectors $b^+$ and $b^-$ into account we arrive at inhomogeneous inequalities, leading to \emph{tropical polyhedra} rather than tropical cones.
This ties in better with ordinary linear optimization; this is the reason for the latter.
Note that $-\infty\1\in\TTmax^n$ trivially satisfies any tropical linear inequality.

The following result shows why tropical linear programming is interesting.
It is concerned with the computational complexity of deciding if some $x\neq-\infty\1$ exists which satisfies all the constraints.
This is the tropical analog to the feasibility problem of linear programming.
\begin{theorem}[{Akian, Gaubert and Guterman \cite{AkianGaubertGutermann12}}]
  Deciding whether a tropical linear program has a nontrivial feasible solution lies in the complexity class $\NP\cap\coNP$.
\end{theorem}
The complexity class $\NP$ comprises those decision problems which admit a nondeterministic polynomial algorithm which can settle the issue in the affirmative.
Its \enquote{mirror image} is the complexity class $\coNP$ which is about nondeterministic polynomial algorithms to certify the negation.
Most crucially, there is no deterministic polynomial algorithm known for tropical linear programming feasibility.
It is a major open problem in theoretical computer science whether or not the complexity classes $\NP\cap\coNP$ and P (deterministic polynomial algorithms) agree.

Yet there is another aspect.
By Theorem~\ref{thm:develin-yu} we may view tropical cones and tropical polytopes as projections of ordinary cones and polytopes over real Puiseux series under the (dual) valuation map.
In other words, each tropical polytope or tropical polyhedron admits a \emph{lift} to their ordinary counterparts, which are defined over real Puiseux series.
The same holds true for tropical linear programs: tropical linear programs arise as projections of ordinary linear programs over Puiseux series.

\begin{example}\label{exmp:long-and-winding}
  For a parameter $r\geq 1$ we define a linear program over $\puiseuxseries{\RR}{t}^*$ with $2r$ variables and $3r+1$ constraints:
  \begin{equation}\label{eq:long-and-winding}
    \begin{array}{l@{\quad}l}
      \text{minimize} & x_1 \\[\jot]
      \text{subject to}
                      & x_1 \leq t^2 \\[\jot] 
                      & x_2 \leq t \\[\jot]
                      & x_{2j+1} \leq t \, x_{2j-1} \, , \quad x_{2j+1} \leq t \, x_{2j} \tikzmark{} \\[\jot]
                      & x_{2j+2} \leq t^{1-1/2^j} (x_{2j-1} + x_{2j}) \tikzmark{} \\[\jot]
                      & x_{2r-1} \geq 0 \, , \quad x_{2r} \geq 0 \enspace .
    \end{array}
    \insertbigbrace{$1 \leq j < r$}
  \end{equation}
  Since all the terms occurring in \eqref{eq:long-and-winding} are nonnegative, the tropicalization is a merely syntactic procedure, term by term.
  This yields:
  \begin{equation}\label{eq:tlaw}
    \begin{array}{l@{\quad}l}
      \text{minimize} & X_1 \\[\jot]
      \text{subject to}
                      & X_1 \leq 2 \\[\jot] 
                      & X_2 \leq 1 \\[\jot]
                      & X_{2j+1} \leq 1 + X_{2j-1} \, , \quad X_{2j+1} \leq 1 + X_{2j} \tikzmark{} \\[\jot]
                      & X_{2j+2} \leq 1-1/2^j + \max (X_{2j-1} , X_{2j}) \tikzmark{} \enspace .
    \end{array}
    \insertbigbrace{$1 \leq j < r$}
  \end{equation}
  We write capital letters for the variables now to mark the distinction from the ordinary linear program \eqref{eq:long-and-winding}.
  The \enquote{built-in positivity} of the tropical semiring was mentioned before; and this is responsible for not needing nonnegativity constraints in \eqref{eq:tlaw}.
\end{example}

While we started out by viewing \eqref{eq:long-and-winding} as a linear program over Puiseux series, we may also consider $t$ a suitable (positive) real number.
In this way, we obtain a family of ordinary linear programs, over the reals, which depend on two parameters, $r$ and $t$.
However, $r$ and $t$ play very different roles.
The idea is that the limit for $t$ approaching infinity is defined, and it is very close to the tropicalization, e.g., measured in the Hausdorff metric.
The precise relationship between tropicalization and this limiting process is known as \emph{Maslov dequantization} \cite{LitvinovMaslovShpiz01}.

As its key property the tropical linear program \eqref{eq:tlp} exhibits convoluted combinatorics.
While the details are beyond the scope of this survey, we wish to give a sketch.
Let $P$ denote the tropical polyhedron which forms the feasible region of \eqref{eq:tlp}.
For a fixed real number $\lambda$ we can consider the \emph{tropical sublevel set} $S_\lambda:=\smallSetOf{x\in P}{x_1\leq \lambda}$.
The latter set is tropically convex, whence, e.g., $u,v\in S_\lambda$ implies $u\maxodot v\in S_\lambda$.
It follows that the coefficientwise maximum in the set $S_\lambda$ is well-defined and unique; let us call this point $p(\lambda)$.
The map $\lambda\mapsto p(\lambda)$ is piecewise linear, and its image is the \emph{tropical central path} of \eqref{eq:tlp}.
One way to express the \enquote{convoluted combinatorics} is to say that, for a certain discrete curvature measure, the tropical central path has very large total curvature.

A standard class of algorithms for solving ordinary linear programs employs the \emph{interior point method}.
It starts at a certain point in the interior of the feasible region, called the \emph{analytic center}, and it follows a branch of a real algebraic curve, known as the \emph{central path} to approach the optimum.
The interior point method is an iterative procedure which employs auxiliary optimization problems at each step; these auxiliary problems are then solved via Newton's method.
There is a host of choices for everything, the precise form of the auxiliary problems, the step sizes and so forth; see Renegar \cite{Renegar:2001} for the details.
In \cite[Theorem 15]{ABGJ:2018} it is shown that for a fairly general class of interior point algorithms the central path (of any linear program over dual real Puiseux series) is mapped to the tropical central path of the tropicalization.
This can be exploited to derive information about the computational complexity of interior point algorithms.
To obtain a precise statement, in the following \emph{log-barrier interior point method} is short for \enquote{primal-dual path-following interior point algorithm with a log-barrier function which iterates in the wide neighborhood of the central path}.
The next result is \cite[Theorem~B]{ABGJ:2018}; see also \cite{ABGJ:SIREV}.
\begin{theorem}\label{thm:long-and-winding}
  The number of iterations required by any log-barrier interior point method, from the analytic center to some point of duality measure $\leq 1$, is exponential in $r$ on the family of linear programs \eqref{eq:long-and-winding} over the reals, provided that $t \gg 1$ is sufficiently large.
\end{theorem}
This theorem is relevant because it is known from the work of Karmarkar \cite{Karmarkar:1984} that some interior point algorithms solve linear programs in polynomial time, measured in the size of the input.
And the log-barrier interior point methods form prominent members of that class.
The \emph{size} is a function of the dimensions $m$ and $n$ as well as the coding lengths of the coefficients of the constraints and the objective function.
Now the linear programs \eqref{exmp:long-and-winding} show that, in general, the running time of the log-barrier interior point methods depend on the coefficient size in a crucial way.
Recently, Theorem~\ref{thm:long-and-winding} was generalized to other barrier functions by Allamigeon, Aznag, Gaubert and Hamdi \cite{AllamigeonAznagGaubertHamdi:2010.10205}.

\subsection{Further applications}
To show how far the ideas of tropical convexity can be carried we briefly mention further connections, without any attempt at completeness.
In particuar, there are even more applications of tropical geometry, beyond tropical convexity, which we ignore here.

\subsubsection*{Matroids}
In Section~\ref{sec:tlinear} we saw how tropical linear spaces arise as tropical convex hulls in some tropical projective space.
When each coefficient of each cocircuit generator is either $0$ or $\infty$, we are in the constant coefficient case; see \cite[\S1.5]{ETC}.
Then the intersection of the tropical linear space with the tropical projective torus is an ordinary polyhedral fan; in this case it is the \emph{Bergman fan} of some matroid; see Feichtner--Sturmfels \cite{FeichtnerSturmfels05}, \cite[\S4.2]{Tropical+Book} and \cite[\S10.8]{ETC}.
Bergman fans occur in the work of Bergman \cite{Bergman:1971} as logarithmic limit sets of algebraic varieties, and they are related to the De Concini--Procesi compactification of the complement of a hyperplane arrangement \cite{DeConciniProcesi:1995}.
More recently, Bergman fans played a central role in groundbreaking work of Adiprasito, Huh and Katz \cite{AHK:2018}, who solved several previously open problems in matroid theory.

\subsubsection*{Bruhat--Tits Buildings}
Keel and Tevelev investigated compactifications of certain moduli spaces in algebraic geometry \cite{KeelTevelev06}.
Fixing some field, $\KK$, equipped with a discrete valuation, this lead them to studying equivalence classes of $R$-lattices, where $R$ is the valuation ring.
The classical setting are the Bruhat--Tits buildings arising from the special linear groups $\SL_d\KK$; see \cite{AbramenkoBrown:2008} for an introduction to buildings.
As a key new tool, Keel and Tevelev introduced \emph{membranes}, which are certain sets of equivalence classes of $R$-lattices; see \cite[\S4]{KeelTevelev06}.
The membranes turn out to agree with the sets of lattice points in certain tropical linear spaces; see \cite{JoswigSturmfelsYu07} and \cite[\S10.9]{ETC}.
More recently, Zhang developed algorithms for computing membranes \cite{Zhang:2021}.
El Maazouz, Hahn, Nebe, Stanojkovski and Sturmfels built on these ideas to study matrix algebras arising in representation theory \cite{MaazouzHahnNebeStanojkovskiSturmfels:2107.00503}.

\subsubsection*{Economics}
Questions in economics which are amenable to methods of tropical convexity are often inherently game-theoretic.
For details see the articles \cite{CohenGaubertQuadrat04,AkianGaubertGutermann12,AkianGaubertVannucci:2108.07748,LohoVegh:2020}, which were mentioned before, and their references.
This line of resaech lead Shiozawa to interprete the classical Ricardian theory of trade \cite{Shiozawa:2015} in terms of tropical convexity; see also \cite[\S9.5]{ETC}.
Shortly after Crowell and Tran characterized incentive compatible mechanisms \cite{CrowellTran:1606.04880}.
Mechanism design is a branch of optimization bordering economics.
The goal is to define auctions for selling items such that, e.g., the revenue of the seller is maximized.

\subsubsection*{Data analysis}
A basic method for dimension reduction in data science is known as principal component analysis (PCA).
Geometrically, a version of this amounts to projecting data points onto a (low-dimensional) affine linear subspace.
Yoshida, Zhang and Zhang proposed a tropical analog to PCA, i.e., a method to project onto a tropical linear space, and applied it to topics in phylogenetic analysis \cite{Yoshida:TPCA}; see also Semple--Steel\cite{SempleSteel:2003} and \cite[\S10.9]{ETC}.
A key goal in phylogenetics is to aggregate data from several trees (which are considered as distinct approximations of one ground truth) into one tree, with some consistency requirements.
This is the \emph{consensus tree problem}; see Bryant \cite{Bryant:2003} and Bryant--Francis--Steel \cite{Bryant+Francis+Steel:2017}.
To this end Lin--Sturmfels--Tang--Yoshida \cite{LinSturmfelsTangYoshida:2017} and Lin--Yoshida \cite{Yoshida:Fermat-Weber} studied a tropical version of the Fermat--Weber problem.
Recently, we found an asymmetric variant which gives stronger results \cite{ComaneciJoswig:2024}.
The key idea of this article is to translate the consensus tree problem into a tropical convex hull problem, where an optimal solution can be found from one (polytropal) cell of the covector decomposition from Theorem \ref{thm:structure}.
This was further generalized by Cox and Curiel \cite{CoxCuriel:2310.07732}.

Zhang, Naitzat and Lim initiated the study of feedforward neural networks with ReLU activation via tropical geometry \cite{ZhangNaitzatLim:1805.07091}.
Later sharp complexity bounds were found by Mont\'ufar--Ren--Zhang \cite{MontufarRenZhang:2104.08135} and Haase--Hertrich--Loho \cite{haase2023lower}.

\subsubsection*{Computational complexity and combinatorial optimization}

Akian, Gaubert and Guterman~\cite{AkianGaubertGutermann12} showed that deciding the feasibility of a tropical linear program is equivalent to computing an optimal strategy in a mean payoff game; see also \cite[Chapter~9]{ETC}.
The latter problem is of interest to researchers in theoretical computer science because it is unknown whether or not such a strategy can be found in polynomial time.
This line of research led Loho to study combinatorial generalizations of tropical linear programming \cite{Loho:2020}.

There are many connections between tropical convexity and combinatorial optimization; the shortest paths and minimum cycle means were already mentioned in Section~\ref{sec:eigen}.
One more recent addition are results concerning the computation of periodic time tables.
This is a notoriously difficult optimization problem with numerous applications, notably in public transport.
Bortoletto, Lindner and Masing developed new methods based on tropical convexity; their article \cite{bortoletto_et_al:OASIcs.ATMOS.2022.3} won a best paper award at the \enquote{22nd Symposium on Algorithmic Approaches for Transportation Modelling, Optimization, and Systems (ATMOS 2022)}.

\bibliographystyle{amsplain}
\bibliography{etc.bib}

\end{document}